\documentclass{birkjour}
 \newtheorem{theorem}{Theorem}[section]
 
 \newtheorem{proposition}[theorem]{Proposition}
 \newtheorem{lemma}[theorem]{Lemma}
 \newtheorem{example}[theorem]{Example}
 \newtheorem{fact}[theorem]{Fact}
 \theoremstyle{definition}
 \newtheorem{definition}[theorem]{Definition}
 \theoremstyle{remark}

 \numberwithin{equation}{section}

\begin{document}

%
%
%
%
%
%
%
%
%

 \title[On proximal uniform normal structure and relatively
nonexpansive mappings]{On proximal uniform normal structure and relatively
nonexpansive mappings}

 \author[Abhik Digar]{Abhik Digar}

 \address{%
 Department of Mathemtics\\
 IIT Ropar\\
 Rupnagar - 140 001\\
 Punjab, India.}

\email{abhikdigar@gmail.com}

 \author[G S Raju Kosuru]{G. Sankara Raju Kosuru}
 \address{%
 Department of Mathemtics\\
 IIT Ropar\\
 Rupnagar - 140 001\\
 Punjab, India.}
 \email{raju@iitrpr.ac.in}





\begin{abstract}
 The main objective of this article is to provide an alternative approach to the central result of [Eldred, A. Anthony, Kirk, W. A., Veeramani, P., Proximal normal structure and relatively nonexpansive mappings, {\bf{Studia Math.}}, vol 171(3), (2005) 283--293] using proximal uniform normal structure. Also
  we provide characterizations of a strictly convex space. 
   Finally, sufficient conditions for the existence of a non-empty proximal sub-pair for a pair in a Banach space are discussed.
\end{abstract}



\maketitle
\section{Introduction and Preliminaries}\label{Section1}
%
%
Let $X$ be a Banach space and let $(A,B)$ be a pair of non-empty subsets of $X.$ We first recall a few terminologies from the literature. The sets $A_0=\{x\in A: \|x-y\|=d(A,B)~\mbox{for some}~y\in B\}$ and $B_0=\{v\in B: \|u-v\|=d(A,B)~\mbox{for some}~u\in A\},$ here $d(A,B)=\inf\{\|x-y\|: x\in A, y\in B\}.$ Then $(A,B)$ is said to be proximal if $A=A_0$ and $B=B_0.$ Also, $(A,B)$ is said to be semisharp proximal (\cite{Raju2011}) if for any $x$ in $A$ (similarly, in $B$), there exists $y, z$ in $B$ (similarly, in $A$) such that $\|x-y\|=d(A,B)=\|x-z\|,$ then it is the case that $y=z.$ In this case, we deonte $y$ by $x'$ (as well as $x$ by $y'$). A pair is known as sharp proximal if it is proximal and semisharp. Further, the pair $(A,B)$ is parallel (\cite{Espinola2008}) if it is sharp proximal and $B=A+h$ for some $h\in X.$ In this case, we have $x'=x+h,~y'=y-h$ for $x\in A, y\in B$ and $\|h\|=d(A,B).$
It is well-known that every non-empty closed convex pair in a strictly convex Banach space is semisharp proximal. Further, every sharp proximal pair in a strictly convex Banach space is proximal parallel.
 We say $(A,B)$ satisfies a certain property if both $A$ and $B$ satisfy the same. If $(A,B)$ is weakly compact, then it is known (\cite{Reich2003}) that $A_0\neq \emptyset,~B_0\neq \emptyset.$ The pair $(A,B)$ is said to have proximal normal structure (\cite{Eldred2005}) if $(H,K)$ is any non-empty closed bounded convex proximal pair in  $(A,B)$ with $\delta (H,K)> d(H,K)$ and $d(H,K)=d(A,B),$ then it is the case that $\delta(p,K)<\delta (H,K)$ and $\delta (q,H)<\delta (H,K)$ for some $(p,q)\in (H,K).$ Here $\delta (H,K)=\sup\{\|x-y\|:x\in A, y\in B\}$ and $\delta (x,K)= \sup\{\|x-y\|: y\in K\}$ for any $x\in H.$ 
A map $T: A\cup B\to A\cup B$ is said to be cyclic if $T(A)\subseteq B,~T(B)\subseteq A.$ Further $T$ is said to be relatively nonexpansive if it
 satisfies $\|Tx-Ty\|\leq \|x-y\|$ for all $x\in A, y\in B.$
A point $x\in A\cup B$ is said to be a best proximity point of a cyclic map $T$ if $\|x-Tx\|=d(A,B)$ and a pair $(u,v)\in (A,B)$ is a best proximity pair of $T$ if $\|u-Tu\|=d(A,B)=\|v-Tv\|.$ Also, a sub-pair $(C,D)$ of $(A,B)$ is said to be $T$-invariant if $T(C)\subseteq D$ and $T(D)\subseteq C.$
 The main result of (\cite{Eldred2005}) is the following:
\begin{theorem}(\cite{Eldred2005})\label{Thm:Eldred2005}
Let $(A,B)$ be a non-empty weakly compact convex pair in a Banach space $X$ and suppose that $(A,B)$ has proximal normal structure. Then every relatively nonexpansive mapping $T$ on $A\cup B$ has a best proximity pair.
\end{theorem}
It has to be observed that the non-emptiness of $(A_0, B_0)$ and the existence of a minimal $T$-invariant subset of $(A_0, B_0)$ play a vital role in proving Theorem \ref{Thm:Eldred2005}. 
In \cite{Eldred2005} the authors used weak compactness to show $A_0\neq \emptyset,~B_0\neq \emptyset$ whereas weak compactness and proximal normal structure for the existence of a minimal $T$-invariant pair.  
 It is natural to ask whether weak compactness is necessary in Theorem \ref{Thm:Eldred2005}. In this regard, we use a geometric notion, called proximal uniform normal structure which has been recently introduced by the authors (\cite{AbhikRaju2021}), to assert that weak compactness can be relaxed.
%

The article is organized as follows. 
 Section \ref{Section 2} deals with characterizations of strict convexity. 
  In Section \ref{Section 3}, we present an alternative approach of Theorem \ref{Thm:Eldred2005}. The final section deals with sufficient conditions for $A_0\neq \emptyset, B_0\neq \emptyset$ for a given non-empty closed bounded convex pair $(A,B).$
  

\section{Characterizations of strictly convex spaces}\label{Section 2}
Let $(A,B)$ be a non-empty pair in  a Banach space $X.$ For $y\in X,$ a sequence $\{x_n\}$ in $A$ is said to be minimizing with respect to $y$ if $\displaystyle\lim_{n\to \infty}\|x_n-y\|=d(y,A).$ Also, a sequence $\{z_n\}$ in $A$ is said to be minimizing with respect to $B$ if $d(z_n, B)\to d(A, B).$ 
 It is well-known that $A$ is approximatively compact if every minimizing sequence in $A$ has a convergent subsequence (\cite{book:Megginson}). 
 We say that $A$ is said to be weak approximatively compact with respect to $B$ if every minimizing sequence in $A$ with respect to $B$ has a weak convergent subsequence in $A.$ The pair $(A,B)$ is weak approximatively compact if $A$ and $B$ are weak approximatively compact with respect to $B$ and $A$ respectively.
The pair $(A,B)$ is said to have the property UC (\cite{Suzuki2009}) if $\{x_n\}, \{y_n\}$ are sequences in $A$ and $\{z_n\}$ is a sequence in $B$ with  $\displaystyle \lim_{n\to \infty} \|x_n-z_n\|=d(A,B)=\displaystyle \lim_{n\to \infty} \|y_n-z_n\|,$ then it is the case that $\displaystyle\lim_{n\to \infty} \|x_n- y_n\|=0.$
 It is easy to see that if a non-empty pair has property UC, then it is semisharp proximal.
The Banach space $X$ is said to have the Kadec-Klee property (\cite{book:Megginson}) if for any sequence on the unit sphere the weak and strong convergence concide. 
The following result provides the characterizations of a strictly convex Banach space.
\begin{theorem}\label{Thm2.1:StrictConvexityEquivalence}
Let $X$ be a Banach space. Then the following statements are equivalent:
\begin{enumerate}
\item[(a)] Every non-empty closed bounded convex pair in $X$ is semisharp proximal.
\item[(b)] $X$ is strictly convex.

\item[(c)] Suppose $(A,B)$ is a non-empty closed bounded convex pair in  $X,~(x,y)\in (A,B)$ and $\{x_n\}$ is a sequence in $A$ such that $\displaystyle\lim_{n\to \infty}\|x_n-y\|=d(A,B)=\|x-y\|.$ Then $\displaystyle\lim_{n\to \infty}x_n=x.$
\item[(d)] If $A$ is a non-empty compact convex subset and $B$ is a non-empty compact subset of $X,$ then $(A,B)$ has property UC.
\item[(e)] If $X$ has the Kadec-Klee property and $(A,B)$ is weak approximatively compact, then $(A,B)$ has property UC.
\end{enumerate}
\end{theorem}
 {\it{Proof.}} We prove that each of statements from $(b)$ to $(e)$ is equivalent to $(a).$ For this, it is enough to prove the following implications: $(a) \Rightarrow (b),~ (a) \Rightarrow (d)$ and $(a) \Leftrightarrow (e).$ The proof of $(b) \Rightarrow (a)$ is straightforward, $(a) \Rightarrow (c)$ follows from Lemma 2.1 of \cite{AbhikRaju2021} and $(c) \Rightarrow (a),~(d) \Rightarrow (a)$ are obvious. 
 
\noindent Suppose $X$ is non-strictly convex but satisfies $(a).$ Then there exists a nontrivial line segment $B=\{tc_1+(1-t)c_2: 0\leq t\leq 1\}$ in $S_X,$ where $c_1, c_2\in S_X$ with $c_1\neq c_2.$ Set $A=\{0\}.$ Then $d(A,B)=1.$ Clearly $B$ is non-empty closed convex subset of $X$ such that $\left\|0-\left(tc_1+(1-t)c_2\right)\right\|=1$ for all $0\leq t\leq 1.$ This contradicts the semisharp proximality of $(A,B).$ Thus we have $(a) \Rightarrow (b).$
%

\noindent For $(a) \Rightarrow  (d),$ let us assume $(a).$ Suppose $A$ is compact convex and $B$ is a compact subset of $X.$ Suppose $\{x_n\}, \{y_n\}$ in $A$ and $\{z_n\}$ in $B$ are sequences such that $\|x_n-z_n\|\to d(A,B)$ and $\|y_n-z_n\|\to d(A,B).$ Since $A,B$ are compact, there exists a subsequence $\{n_k\}$ of $\{n\}$ such that $x_{n_k}\to x_0, y_{n_k}\to y_0$ and $z_{n_k}\to z_0$ for some $x_0, y_0$ in $A$ and $z_0$ in $B.$ Then it is clear that $\|x_0-z_0\|=d(A,B)=\|y_0-z_0\|.$ By semisharp proximality, $x_0=y_0.$ This shows that $(A,B)$ has property UC. 

\noindent $(a) \Leftrightarrow  (e):$ Assume $(a).$ Let $(A,B)$ be a bounded convex semisharp proximal pair in $X.$ Suppose $(A,B)$ is weak approximatetively compact, and $X$ has the Kadec-Klee property. Let $(x_n),(y_n)$ in $A$ and $(z_n)$ in $B$ be sequences such that $\|x_n- z_n\|\to d(A,B)$ and $\|y_n-z_n\|\to d(A,B).$ It is clear that $(x_n), (y_n)$ in $A$ and $(z_n)$ in $B$ are minimizing sequences.
Then there exists a subsequence $(n_k)$ of $(n)$ such that $x_{n_k}\to x, y_{n_k}\to y$ and $z_{n_k}\to z$ weakly for some $x, y\in A$ and $z\in B.$ Then $x_{n_k}-z_{n_k}\to x-z$ and $y_{n_k}-z_{n_k}\to y-z$ weakly. Since the norm is weakly lower semicontinuous, we have $\|x_{n_k}-z_{n_k}\|\to \|x-z\|$ and $\|y_{n_k}-z_{n_k}\|\to \|y-z\|.$ By the Kadec-Klee property of $X,$ it is clear that $x_{n_k}-z_{n_k}\to x-z$ and $y_{n_k}-z_{n_k}\to y-z$ strongly. Thus $x_{n_k}-y_{n_k}\to x-y.$ Moreover, $\|x-z\|=d(A,B)=\|y-z\|.$ By $(a), x=y.$ Hence $\|x_{n_k}-y_{n_k}\|\to 0.$ 

\noindent To see $(e) \Rightarrow  (a),$ suppose $x,y\in A$ and $z\in B$ such that $\|x-z\|=d(A,B)=\|y-z\|.$ Set $A=\overline{co}(\{x,y\})$ and $B=\{z\}.$ Then $(A,B)$ is a closed bounded weak approximatively compact convex pair. 
 By property UC, we have $x=y.$ \qed
%
%
%
%
%
%

We recall proximal uniform normal structure from \cite{AbhikRaju2021}.
For a non-empty pair $(A,B)$ in $X,$ define a collection $\Upsilon (A,B)$ consisting of non-empty closed bounded convex proximal subsets $(H_1, H_2)$ of $(A,B)$ with at least one of $H_1$ and $H_2$ has more than one point and $d(H_1, H_2)=d(A,B).$ If $r(H_1, H_2)=\inf\{\delta (x,H_2): x\in H_1\},~r(H_2, H_1)=\inf\{\delta (y,H_1): y\in H_2\},$ then define $R(H_1,H_2)=\max\{r(H_1, H_2), r(H_2, H_1)\}.$ Also, $r(H_1)=r(H_1,H_1)$ and $\delta(H_1)=\delta(H_1,H_1).$ Define $N(A,B)=\displaystyle\sup\left\{\frac{R(H_1, H_2)}{\delta(H_1, H_2)}:(H_1, H_2)\in\Upsilon (A,B) \right\}.$
 \begin{definition} \label{Def:PUNS}
 A non-empty bounded convex pair $(A,B)$ is said to have proximal uniform normal structure if $N(A,B)<1.$
 \end{definition} 
%
We say that $X$ has proximal uniform normal structure if there exists $c<1$ such that $N(A,B)\leq c$ for every non-empty closed bounded convex proximal pair $(A,B)$ of $X$ with one of $A,B$ contains more than one element.
It is known (\cite{AbhikRaju2021}) that if a non-empty closed bounded convex pair in a Banach space has proximal uniform normal structure, then it is semisharp proximal. Thus 
in view of Theorem \ref{Thm2.1:StrictConvexityEquivalence}, we have the following theorem. 
 \begin{theorem}\label{Thm2.2:Thm:StrictConvexity}
 Let $X$ be a Banach space having proximal uniform normal structure. Then $X$ is strictly convex.
 \end{theorem}
It has to be observed that Theorem \ref{Thm2.1:StrictConvexityEquivalence} discusses on geometrical structures of a strictly convex Banach space. In this view point, it is interesting to study the converse of Theorem \ref{Thm2.2:Thm:StrictConvexity}. But the same is not known even in a uniformly convex Banach space setting (see Conjecture 5.3 in \cite{AbhikRaju2021}). However, we prove this in the setting of a Hilbert space. A semisharp proximal pair $(A,B)$ in a Banach space is said to have the Pythagorean property (\cite{EspinolaRaju2015}) if for each $(x,y)\in (A_0,B_0),$ we have $\|x-y'\|^2+\|x-x'\|^2=\|x-y\|^2=\|x'-y\|^2+\|y'-y\|^2.$ It is well-known that every non-empty closed convex pair in a Hilbert space has the Pythagorean property.
%
\begin{theorem}\label{Thm2.3:HilbertSpace}
Every non-empty closed bounded convex proximal pair in a Hilbert space has proximal uniform normal structure.
\end{theorem} 
{\it{Proof.}} Let $(A,B)$ be a non-empty closed bounded convex proximal pair in a Hilbert space $\mathcal{H}.$ Suppose $\Upsilon(A,B)$ be the collection of all closed bounded convex proximal pair $(K_1, K_2)$ of subsets of $(A,B)$ such that one of $K_1, K_2$ contains more than one element and $d(K_1, K_2)=d(A,B).$ 
 As every $(K_1, K_2)$ in $\Upsilon(A,B)$ is proximal parallel, we have $r(K_1)=r(K_2)$ and $\delta(K_1)=\delta(K_2).$ Fix $(K_1, K_2)\in \Upsilon(A,B)$ and let $x\in K_1.$ Then for any $y\in K_2,$ by Pythagorean property we have $\|x-y\|^2=\|x-y'\|^2+d^2(A,B).$ This proves that $\delta^2(x,K_2)=\delta^2(x,K_1)+d^2(A,B).$ 
Therefore, $r^2(K_1,K_2)= r^2(K_1)+d^2(A,B)$ and $\delta^2(K_1, K_2)=\delta^2(K_1)+d^2(A,B).$ Similarly, $r^2(K_2, K_1)= r^2(K_2)+d^2(A,B)=r^2(K_1)+d^2(A,B).$ As every Hilbert space has uniform normal structure (\cite{book:Goebel}), there exists $c\in (0,1)$ (independent of $(K_1, K_2)$) such that $r(K_1)\leq c \delta(K_1).$ Therefore, 
\begin{eqnarray*}
 \displaystyle \frac{R(K_1, K_2)}{\delta(K_1, K_2)}= \frac{\sqrt{r^2(K_1)+d^2(A,B)}}{\sqrt{\delta^2(K_1)+d^2(A,B)}}\leq \frac{\sqrt{c\delta^2(K_1)+d^2(A,B)}}{\sqrt{\delta^2(K_1)+d^2(A,B)}}.
\end{eqnarray*}
This proves that $N(A,B)<1.$ \qed

\section{Existence of Best Proximity Pairs}\label{Section 3}

In this section we prove the existence of a best proximity pair for a relatively nonexpansive mapping. In literature, such existence theorems are discussed using weak compactness and proximal normal structure (Theorem 2.1 of \cite{Eldred2005}, Theorem 3.11 of \cite{Espinola2008}, Theorem 3.12 of \cite{Raju2010} etc.).
Here we replace weakly compact sets by closed bounded sets and proximal normal structure by proximal uniform normal structure.
The following technical results are useful for our main result.
  \begin{proposition}\label{Prop3.1:NondiametralProximal}
 Let $(A,B)$ be a non-empty closed bounded convex proximal pair in a Banach space $X.$ Suppose $(A,B)$ has proximal uniform normal structure. Then there exists $c\in (0,1)$ 
with the following property: 

for every $(H_1, K_1)\in \Upsilon(A,B)$ there exists $(z_1, z_2)\in (H_1,K_1)$ with $\|z_1-z_2\|=d(H_1,K_1)$ such that $\max\{\delta(z_1, K_1), \delta(z_2, H_1)\}\leq c \delta(H_1,K_1).$
 \end{proposition}
{\it{Proof.}}
 Set $\epsilon=\frac{1-N(A,B)}{2}.$ 
 Then for $(H_1, K_1)\in \Upsilon(A,B),$ there exists a pair $(u,v)\in (H_1, K_1)$ such that 
\begin{eqnarray*}
\max \{\delta(u, K_1), \delta(v, H_1)\}\leq (N(A,B)+\epsilon) \delta(H_1, K_1)
\end{eqnarray*}
By setting $z_1=\frac{u+v'}{2}$ and $z_2=\frac{u'+v}{2},$ we have $\|z_1-z_2\|=d(A,B).$
 Now, for any $y\in K_1,$ 
\begin{eqnarray*}
\|z_1-y\| &= & \left\|\frac{u+v'}{2}-y \right\|\\
&\leq & \left\| \frac{u-y}{2}\right\|+\left\| \frac{v'-y}{2}\right\|\\
&\leq & \frac{1}{2} (N(A,B)+\epsilon)\delta(H_1, K_1)+\frac{1}{2}\delta(H_1, K_1)\\
&=&  c \delta(H_1, K_1).
\end{eqnarray*} 
%
\noindent  Here $c=\frac{N(A,B)+\epsilon+1}{2}=\frac{3+N(A,B)}{4}<1.$ Therefore, $\delta(z_1,K_1)\leq c \delta(H_1, K_1).$ Analogously, $\delta(z_2,H_1)\leq c \delta(H_1, K_1).$ This completes the proof. \qed

  \begin{lemma}\label{Lem3.1:MainLemma}
 Let $(A,B)$ be as in Proposition \ref{Prop3.1:NondiametralProximal}.
  Suppose $T$ is a relatively nonexpansive mapping on $A\cup B.$ Then there exist a $T$-invariant pair $(H_1, K_1)$ in $\Upsilon(A,B)$ and a constant $c\in (0,1),$ independent of $(H_1, K_1)$  such that $\delta(H_1, K_1)-d(A,B)\leq c\displaystyle \left(\delta(A,B)-d(A,B)\right).$ 
 \end{lemma}
{\it{Proof.}}
 By Proposition \ref{Prop3.1:NondiametralProximal}, there exists $c\in (0,1)$
  and a pair $(z_1, z_2)\in (A, B)$ such that $\|z_1-z_2\|=d(A,B)$,
 \begin{align}\label{Eqn3.1}
 \begin{split}
 \delta(z_1, B)-d(A,B) &\leq  c \left(\delta(A,B)-d(A,B)\right),\\
 \delta(z_2, A)-d(A,B) &\leq  c \left(\delta(A,B)-d(A,B)\right).
 \end{split}
 \end{align}
 \noindent  Let $\Omega$ be the collection of $L_\alpha\subseteq A\cup B,~\alpha \in \Lambda ,$ for some index set $\Lambda$ with $(L_\alpha\cap A, L_\alpha\cap B)$ is a non-empty closed convex $T$-invariant proximal pair in $(A,B)$ 
  such that $ (z_1, z_2)\in (L_\alpha\cap A, L_\alpha\cap B).$ 
\noindent  Then $\Omega \neq \emptyset,$ since $A\cup B\in \Omega.$ Set $L=\bigcap \{L_\alpha: L_\alpha\in \Omega\}.$ Thus 
\begin{eqnarray*}
L\cap A = \bigcap \{L_\alpha\cap A: L_\alpha\in \Omega\},~
L\cap B = \bigcap \{L_\alpha\cap B: L_\alpha\in \Omega\}.
\end{eqnarray*}
 Clearly, $(L\cap A, L\cap B)$ is a closed convex $T$-invariant proximal pair and $(z_1, z_2)\in (L\cap A, L\cap B).$ 
Set
 $A_1= T(L\cap B)\bigcup \{z_1\}; ~B_1= T(L\cap A)\bigcup \{z_2\}.$
Clearly,
\begin{align}\label{subset1}
\overline{co}(A_1)\subseteq L\cap A~~~\mbox{and}~~~\overline{co}(B_1)\subseteq L\cap B.
\end{align}
\noindent Here $\overline{co}(E)$ is the closed convex hull of $E$ for $E\subset X.$ It can be easily verified that $\left(\overline{co}(A_1), \overline{co}(B_1)\right)$ is $T$-invariant.
We now prove that $(\overline{co}(A_1), \overline{co}(B_1))$ is proximal. Let $x\in \overline{co}(A_1).$ Then there exists a sequence $\{x_n\}$ in the convex hull of $A,~co(A_1)$ that converges, say to $x.$ For $n\in \mathbb{N},$ there exist $N\in \mathbb{N},~\lambda_i^{(n)}\geq 0,~ (1\leq i\leq N+1),~\displaystyle \sum_{i=1}^{N+1}\lambda_i^{(n)}=1$ and $y_j^{(n)}\in L\cap B,~(1\leq j\leq N)$ such that $x_n=\displaystyle \sum_{i=1}^{N}\lambda_i^{(n)} T(y_i^{(n)})+\lambda_{N+1}^{(n)}z_1.$ As $(L\cap A, L\cap B)$ is proximal, we get $x_j^{(n)}\in L\cap A$ with $\|x_j^{(n)}-y_j^{(n)}\|=d(A,B)$ for $1\leq j\leq N.$ By setting $ y_n=\displaystyle \sum_{i=1}^{N}\lambda_i^{(n)} T(x_i^{(n)})+\lambda_{N+1}^{(n)}z_2$ we have $\|x_n-y_n\|=d(A,B).$
Thus there exists $y\in L\cap B$ such that $\|x-y\|=d(A,B).$ As $\|x-y_n\|\leq \|x-x_n\|+\|x_n-y_n\|\to d(A,B),$
by Theorem \ref{Thm2.1:StrictConvexityEquivalence},
  the sequence $\{y_n\}$ converges to $y.$
 Consequently, $y\in \overline{co}(B_1).$ Similarly,
   one can show that for a given $v\in \overline{co}(B_1),$ there exists $u\in \overline{co}(A_1)$ such that $\|u-v\|=d(A,B).$ Therefore, $(\overline{co}(A_1), \overline{co}(B_1))$ is a non-empty closed convex proximal pair. 
\noindent  Since $(z_1, z_2)\in (\overline{co}(A_1), \overline{co}(B_1)),$ we have $d(\overline{co}(A_1), \overline{co}(B_1))\leq \|z_1-z_2\|=d(A,B)$ and hence $(\overline{co}(A_1), \overline{co}(B_1))\in \Omega.$ Thus
\begin{align}\label{subset2}
L\cap A\subseteq \overline{co}(A_1),~~ L\cap B\subseteq \overline{co}(B_1).
\end{align}
Combining (\ref{subset1}) and (\ref{subset2}), we have
$L\cap A= \overline{co}(A_1),~~~ L\cap B= \overline{co}(B_1).$
Define
{\scriptsize{
 \begin{align*}
H_1=\{x\in L\cap A: \max\{\delta(x, L\cap B)-d(A,B), \delta(x', L\cap A)-d(A,B)\}\leq c (\delta(A,B)-d(A,B))\},\\
K_1=\{y\in L\cap B: \max\{\delta(y, L\cap A)-d(A,B),\delta(y', L\cap B)-d(A,B)\}\leq c (\delta(A,B)-d(A,B))\}.
 \end{align*}}}
By (\ref{Eqn3.1}), we have $(z_1, z_2)\in (H_1, K_1).$ It is easy to see that $(H_1, K_1)$ is convex and proximal. By Theorem \ref{Thm2.1:StrictConvexityEquivalence}, one can show that $(H_1, K_1)$ is closed.
 Now we prove that $(H_1, K_1)$ is $T$-invariant. Let $w\in H_1.$ 
  Suppose $u\in L\cap A$ and $\epsilon>0.$ Since $\overline{co}(A_1)=L\cap A,$ there exist $N\in \mathbb{N},~\alpha_i\geq 0~(\mbox{for}~1\leq i\leq N+1),~\displaystyle \sum_{i=1}^{N+1}\alpha_i=1$ and $y_j\in L\cap B~(\mbox{for}~1\leq j\leq N)$ such that $\displaystyle\left\|\sum_{i=1}^{N} \alpha_i T(y_i)+\alpha_{N+1}z_1-u\right\|<\epsilon.$ 
Then we have
{\scriptsize{
\begin{eqnarray*}
\|Tw-u\|-d(A,B)&\leq &\left\|Tw-\displaystyle\sum_{i=1}^{N} \alpha_i T(y_i)-\alpha_{N+1}z_1\right\|-d(A,B)+\epsilon\\
&\leq & \displaystyle\sum_{i=1}^{N} \alpha_i  \left\|Tw-T(y_i)\right\|+\alpha_{N+1}\left\|Tw-z_1\right\|-d(A,B)+\epsilon\\
&\leq & \displaystyle\sum_{i=1}^{N} \alpha_i  \left\|w-y_i\right\|+\alpha_{N+1}\left\|Tw-z_1\right\|-d(A,B)+\epsilon\\
&=& \displaystyle\sum_{i=1}^{N} \alpha_i  \left(\left\|w-y_i\right\|-d(A,B)\right)+\alpha_{N+1}\left(\left\|Tw-z_1\right\|-d(A,B)\right)+\epsilon\\
&\leq &\displaystyle\sum_{i=1}^{N} \alpha_i c (\delta (A,B)-d(A,B))+\alpha_{N+1}c (\delta (A,B)-d(A,B))+\epsilon\\
&=& c (\delta (A,B)-d(A,B))+\epsilon.
\end{eqnarray*}}}
In a similar fashion, 
 we can show that $\|Tw'-u'\|-d(A,B)\leq c (\delta (A,B)-d(A,B))+\epsilon.$
Since $\epsilon >0$ was chosen arbitrarily, $Tw\in K_1.$ Thus $(H_1, K_1)$ is $T$-invariant.
This completes the proof of the lemma. \qed

It has to be observed that every $c \in \displaystyle \left[\frac{3+N(A,B)}{4},1\right)$ satisfies the conclusions of Proposition \ref{Prop3.1:NondiametralProximal} and Lemma \ref{Lem3.1:MainLemma}. However, we fix $c=\frac{3+N(A,B)}{4}$ to prove our main results.
Again applying Proposition \ref{Prop3.1:NondiametralProximal} for the pair $(H_1,K_1),$ we have for every $(H_2,K_2)\in \Upsilon(H_1,K_1)$ there exists $(u_1,u_2)\in (H_2,K_2)$ with $\|u_1-u_2\|=d(H_2,K_2)$ such that $\max\{\delta(u_1, K_2), \delta(u_2, H_2)\}\leq c_1 \delta(H_2,K_2)$ for some $c_1\in (0,1).$ Here $c_1=\frac{3+N(H_1,K_1)}{4}.$ 
As $N(H_1,K_1)\leq N(A,B),$ we may choose $c_1=c.$ Thus, 
for every $(H_1, K_1)\in \Upsilon(A,B)$ and $(H_2,K_2)\in \Upsilon(H_1,K_1)$ there exists $(w_1, w_2)$ in $(H_2,K_2)$ with $\|w_1-w_2\|=d(H_2,K_2)$ such that $\max\{\delta(w_1, K_2), \delta(w_2, H_2)\}\leq c \delta(H_2,K_2).$

 \begin{theorem}\label{Thm3.1:MainTheorem}
 Let $(A,B)$ be a non-empty closed bounded convex proximal pair in  a Banach space having proximal uniform normal structure. 
  Then every relatively nonexpansive mapping $T$ on $A\cup B$ has a best proximity pair. 
  \end{theorem}
{\it{Proof.}}
  By Lemma \ref{Lem3.1:MainLemma}, there exist $c\in (0,1)$ and a 
   pair $(H_1, K_1)\in \Upsilon(A,B)$ such that 
   $ \delta(H_1, K_1)-d(A,B)\leq c(\delta(A,B)-d(A,B)).$ Again by applying Lemma \ref{Lem3.1:MainLemma}, there exists a non-empty closed convex proximal $T$-invariant pair $(H_2, K_2)\subseteq (H_1, K_1)$ such that $d(H_2, K_2)=d(H_1, K_1)=d(A,B)$ and $\delta(H_2, K_2)-d(H_2, K_2) \leq  c(\delta(H_1, K_1)-d(H_1,K_1)).$ Then $\delta(H_2, K_2)-d(A,B)\leq  c^2(\delta(A,B)-d(A,B)).$
\noindent Thus, for $n\geq 1,$ there exists a non-empty closed convex proximal $T$-invariant pair $(H_{n}, K_{n})\subseteq (A,B)$ such that $d(H_{n}, K_{n})=d(A,B)$ and
\begin{eqnarray*}
\delta(H_{n}, K_{n})-d(A,B)\leq c^{n}(\delta(A,B)-d(A,B)).
\end{eqnarray*}
  Hence $\delta(H_{n}, K_{n})\to d(A,B).$
We claim that for some $n\geq 1,~H_n$ is singleton.
  Otherwise, for all $n\geq 1,~(H_n,K_n) \in\Upsilon(A,B)$ with $\delta(H_n,K_n)\to d(A,B).$ Then $R(H_n,K_n)\to d(A,B).$ This contradicts the hypothesis that $N(A,B)<1.$   
Thus there exists $n_0\in \mathbb{N}$ such that $H_{n_0}$ is singleton, say $\{x_0\}$ and hence $K_{n_0}=\{Tx_0\}.$ This completes the proof. \qed

\section{Sufficient Conditions for $A_0\neq \emptyset,~B_0\neq \emptyset$}\label{Section 4}

It is well-known that a non-empty weakly compact (or approximatively compact) pair $(A,B)$ in a Banach space
 has a proximal sub-pair (\cite{Reich2003, Raju2015, book:Singer}).
%
Also, if $(A,B)$ is a non-empty closed bounded pair in a Banach space such that $(A,B)$ is weak approximatively compact, $ A_0\neq \emptyset,~B_0\neq \emptyset.$
It is known that the non-emptiness of $(A_0,B_0)$ is a necessary condition for a cyclic map on $(A,B)$ to have a best proximity point. Here we aim to show that for every non-empty pair $(A,B)$ having proximal uniform normal structure, $(A_0,B_0)$ is non-empty. In fact, we prove a little stronger result. Consider $\overline{\Upsilon}(A, B),$ the collection of pair $(H_1, H_2)$ of all non-empty closed bounded convex subsets of $(A,B)$ such that one of $H_1, H_2$ contains more than one element and $d(H_1, H_2)=d(A,B).$ Suppose
$c_0:=\displaystyle\sup\left\{\frac{R(C_1, C_2)}{\delta(C_1, C_2)}: (C_1, C_2)\in\overline{\Upsilon}(A,B)\right\}.$
%
The following fact is immediate. 
\begin{fact}\label{Fact4.1}
Let $(A,B), \overline{\Upsilon}(A,B)$ and $c_0$ be as above. Suppose $c_0<1.$ Then there exists $c\in (c_0,1)$ such that for $(C_1,C_2)\in \overline{\Upsilon}(A,B),$ there is a pair $(z_1, z_2)\in (C_1,C_2)$ satisfying
$\delta(z_1, C_2)\leq c\delta(C_1, C_2)~~\mbox{and}~~
\delta(z_2, C_1)\leq c\delta(C_1, C_2).$
\end{fact} 
%
Now we prove the main result of this section. 
\begin{theorem}\label{Thm4.1:ExistenceProximal2}
Let $(A,B)$ be a non-empty closed bounded convex pair in a Banach space $X.$ Suppose
$c_0<1,$ where $c_0$ is as above.
 Then $A_0\neq \emptyset,~B_0\neq \emptyset.$  
\end{theorem}  
{\it{Proof.}}
There exists a pair of sequences $\left\{(x_n, y_n)\right\}$ in $(A,B)$ such that $\|x_n-y_n\|\to d(A,B).$ If $\{(x_n, y_n)\}$ is a constant sequence, then the theorem holds tivially. Therefore, without loss of generality, we assume that $\{(x_n, y_n)\}$ is non-constant.
For $n\in \mathbb{N},$ set $H_n^0=\overline{co}\left\{x_n, x_{n+1}, x_{n+2}, \cdots\right\}$ and $K_n^0=\overline{co}\left\{y_n, y_{n+1}, y_{n+2}, \cdots\right\}.$ Then $\left(H_n^0, K_n^0\right)$ is a decreasing sequence in $\overline{\Upsilon}(A,B).$
 By Fact \ref{Fact4.1}, there exists $c\in (c_0,1)$ such that $\mathfrak{A}(C_1)\neq \emptyset, \mathfrak{A}(C_2)\neq \emptyset$ for all $(C_1, C_2)\in \overline{\Upsilon}(A,B).$ Here 
\begin{eqnarray*}
\mathfrak{A}(C_1)&=& \left\{x\in C_1:\delta(x, C_2)-d(A,B)\leq c\left(\delta(C_1, C_2)-d(A,B)\right)\right\};\\
\mathfrak{A}(C_2)&=& \left\{y\in C_2:\delta(y, C_1)-d(A,B)\leq c\left(\delta(C_1, C_2)-d(A,B)\right)\right\}.
\end{eqnarray*}
For $n\in \mathbb{N},$ define
\begin{eqnarray*}
H_n^1= \overline{co}\left(\bigcup_{j=n}^{\infty} \mathfrak{A}(H_j^0)\right)~~\mbox{and}~~K_n^1=\overline{co}\left(\bigcup_{j=n}^{\infty} \mathfrak{A}(K_j^0)\right).
\end{eqnarray*} 
Note that $d\left(\mathfrak{A}(H_n^0),\mathfrak{A}(K_n^0)\right)\to d(A,B)$ and this leads to obtain $d(H_n^1,K_n^1)=d(A,B),~ n\geq 1.$ 
Thus $d(H_n^1, K_n^1)$ is a non-empty closed bounded convex pair of subsets of $(H_n^0, K_n^0)$ with $d(H_n^1, K_n^1)=d(H_n^0, K_n^0)=d(A,B),~n\geq 1.$ If both of $H_n^1, K_n^1$ contain exactly one point for some $n\in \mathbb{N},$ then the theorem holds. So we may assume that at least one of $H_n^1, K_n^1$ contains more than one point for each $n\in \mathbb{N}.$ 
  For a fixed $n\geq 1,$ let $x\in \displaystyle \bigcup_{j=n}^{\infty} \mathfrak{A}(H_j^0)$ and $y\in \displaystyle \bigcup_{j=n}^{\infty} \mathfrak{A}(K_j^0).$ There exist $p, q\geq n$ such that $x\in \mathfrak{A}(H_p^0)$ and $y\in \mathfrak{A}(K_q^0).$ Without loss of generality, we may assume that $p\leq q.$ Clearly,
  $\delta(x, K_p^0)-d(A,B)\leq c(\delta(H_p^0, K_p^0)-d(A,B))$ and $\delta(y, H_q^0)-d(A,B)\leq c(\delta(H_p^0, K_p^0)-d(A,B)).$ It follows that $\delta\left(\displaystyle \bigcup_{j=n}^{\infty}\mathfrak{A}(H_j^0), \displaystyle \bigcup_{j=n}^{\infty}\mathfrak{A}(K_j^0)\right)-d(A,B)\leq c(\delta(H_p^0, K_p^0)-d(A,B))$ and hence $\delta(H_n^1, K_n^1)-d(A,B)\leq c (\delta(H_n^0, K_n^0)-d(A,B)).$ 
%

\noindent Continuing this process, for any $m, n\geq 1,$ we can construct a sequence of non-empty closed bounded convex pair $(H_n^m, K_n^m)$ in $(A,B)$
with $d(H_n^m, K_n^m)=d(A,B).$
 If both of $H_n^m, K_n^m$ contain exactly one point for some $m, n\geq 1,$ then the theorem holds. Suppose at least one of $H_n^m, K_n^m$ contains more than one point for each $m, n\geq 1.$ Then we get
\begin{eqnarray}\label{equation1}
\delta(H_n^m, K_n^m)-d(A,B)\leq c(\delta(H_n^{m-1}, K_n^{m-1})-d(A,B)).
\end{eqnarray}
Also, we have the following inclusions: for $m\geq 0, n\geq 1,$
\begin{eqnarray*}
(H_1^m, K_1^m)&\supseteq &(H_2^m, K_2^m))\supseteq \cdots \supseteq (H_n^m, K_n^m)\supseteq \cdots;\\
(H_n^0, K_n^0)&\supseteq &(H_n^1, K_n^1))\supseteq \cdots \supseteq (H_n^m, K_n^m)\supseteq \cdots.
\end{eqnarray*}
From (\ref{equation1}), we get $\delta(H_n^m, K_n^m)-d(A,B)\leq c^m(\delta(H_1^{0}, K_1^{0})-d(A,B))$
for $m\geq 0,~n\geq 1.$ 
Hence 
 we have $\delta(H_{n+1}^n, K_{n+1}^n)\to d(A,B)$ as $n\to \infty.$ It follows that $R(H_{n+1}^n, K_{n+1}^n)\to d(A,B)$ as $n\to \infty.$ This forces $c_0=1,$ a contradiction. \qed
%
\begin{example}
Consider the sequence space $\ell^2.$
  Suppose $0\neq h\in \mathbb{R}$ and $A=\Big\{(x_1, x_2, x_3,\cdots)\in \ell^2:  \displaystyle \sum_{i=2}^{\infty}|x_i|^2\leq 1, x_1=0\Big\}, B=\Big\{(y_1, y_2, y_3, \cdots)\in \ell^2: \displaystyle \sum_{i=2}^{\infty}|y_i|^2\leq 1, y_1=h\Big\}.$ By Theorem \ref{Thm2.3:HilbertSpace}, $(A,B)$ has proximal uniform normal structure.
  It can be verified that $A_0=A$ and $B_0=B.$
\end{example}
It has to be observed that the condition $c_0<1$ in Theorem \ref{Thm4.1:ExistenceProximal2} is sufficient but not necessary.
\begin{example}
Let $A=\{(0,x)\in \mathbb{R}^2: 0\leq x\leq 1\}$ and $B=\{(1,y)\in \mathbb{R}^2: 0\leq y\leq 1\}$ be two subsets in the Banach space $\left(\mathbb{R}^2, \|\cdot\|_\infty\right).$ An easy calculation shows that $A_0=A$ and $B_0=B.$ Notice that $c_0=1.$
\end{example}

\end{document}